\documentclass[11pt]{amsart}

\usepackage{latexsym} \usepackage{amsmath}
\newtheorem{thm}{Theorem}[section]
\newtheorem{lem}[thm]{Lemma}
\newtheorem{cor}[thm]{Corollary}

\theoremstyle{definition}
\newtheorem{defn}[thm]{Definition}

\newtheorem{rem}[thm]{Remark}

\begin{document}

\title[Acylindrical accessibility for groups acting on $\mathbb R$-trees]
{Acylindrical accessibility for groups acting on $\mathbb R$-trees}

\author[I.~Kapovich]{Ilya Kapovich}

\address{Dept. of Mathematics, University of Illinois at Urbana-Champaign, 1409
West Green Street, Urbana, IL 61801, USA} \email{kapovich@math.uiuc.edu}

\thanks{The first author was supported by
the U.S.-Israel Binational Science Foundation grant BSF-1999298}

\author[R.Weidmann]{Richard Weidmann} \address{ Fachbereich Mathematik,
Johann Wolgang Goethe-Universit{\"a}t,
Robert Mayer-Strasse 6-8,
60325 Frankfurt (Main), Germany} \email{rweidman@math.uni-frankfurt.de}

\subjclass[2000]{20F67}

\keywords{}

\date{\today}

\begin{abstract} We prove an acylindrical accessibility theorem
for finitely generated groups acting on $\mathbb R$-trees. Namely,
we show that if $G$ is a freely indecomposable non-cyclic
$k$-generated group acting minimally and $D$-acylindrically on an
$\mathbb R$-tree $X$ then there is a finite subtree
$T_{\varepsilon}\subseteq X$ of measure at most
$2D(k-1)+\varepsilon$ such that $GT_{\varepsilon}=X$. This
generalizes theorems of Z.~Sela and T.~Delzant about actions on
simplicial trees.
\end{abstract}

\maketitle

\section{Introduction}

An isometric action of a group $G$ on an $\mathbb R$-tree $X$ is
said to be \emph{$D$-acylindrical} (where $D\ge 0$) if for any
$g\in G, g\ne 1$ we have $diam\ {\rm Fix}(g) \le D$, that is any
segment fixed point-wise by $g$ has length at most $D$. For
example the action of an amalgamated free product $G=A \ast_C B$
on the corresponding Bass-Serre tree is $2$-acylindrical if $C$ is
malnormal in $A$  and $1$-acylindrical if $C$ is malnormal in both
$A$ and $B$. In fact the notion of acylindricity  seems to have
first appeared in this context in the work of Karras and
Solitar~\cite{KarS}, who termed it \emph{being $r$-malnormal}.

Sela~\cite{Se1} proved an important acylindrical accessibility
result for finitely generated groups which, when applied to
one-ended groups, can be restated as follows: for any one-ended
finitely generated group $G$ and any $D\ge 0$ there is a constant
$c(G,D)>0$ such that for any minimal $D$-acylindrical action of
$G$ on a simplicial tree $X$ the quotient graph $X/G$ has at most
$c(G,D)$ edges. This fact plays an important role in Sela's theory
of JSJ-decompositions for word-hyperbolic groups~\cite{Se2} and
thus in his solution of the isomorphism problem for torsion-free
hyperbolic groups~\cite{Se}. Moreover, acylindrical splittings
feature prominently in relation to the Combination Theorem of
Bestvina-Feighn~\cite{BF92,BF96} and its various applications and
generalizations~\cite{Da,Gi96,K01,KM98}. Unlike other kinds of
accessibility results, such as Dunwoody
accessibility~\cite{Du85,Du93}, Bestvina-Feighn generalized
accessibility~\cite{BF91,BF91a} and strong accessibility
(introduced by Bowditch~\cite{Bow98a} and proved by Delzant and
Potyagailo~\cite{DP}), acylindrical accessibility holds for
finitely generated and not just finitely presented groups.

Delzant~\cite{D3} obtained a relative version of Sela's theorem
for finitely presented groups with respect to a family of
subgroups. In particular, he showed that the constant $c(G,D)$
above can be chosen to be $12DT$, where $T$ is the number of
relations in any finite presentation of $G$ where all relators
have length three. Weidmann~\cite{W} used the theory of Nielsen
methods for groups acting on simplicial trees to show that for any
$k$-generated one-ended group one can choose $c(G,D)=2D(k-1)$. In
the present paper we obtain an analogue of this last result for
groups acting on $\mathbb R$-trees.

Before formulating our main result let us recall the notion of Nielsen
equivalence:

\begin{defn}[Nielsen equivalence] Let $G$ be a group and let $M=(g_1,\dots
,g_n)\in G^n$ be an $n$-tuple of elements of $G$. The following moves
are called
\emph{elementary Nielsen moves} on $M$:

\begin{enumerate}
\item[(N1)] For some $i, 1\le i\le n$ replace $g_i$ by
$g_i^{-1}$ in $M$.
\item[(N2)] For some $i\ne j$, $1\le i,j\le n$ replace $g_i$
by $g_ig_j$ in $M$.
\item[(N3)] For some $i\ne j$, $1\le i,j\le n$ interchange
$g_i$ and $g_j$ in $M$.
\end{enumerate}

We say that $M=(g_1,\dots, g_n)\in G^n$ and $M'=(f_1,\dots , f_n)\in G^n$ are
\emph{Nielsen-equivalent}, denoted $M\sim_N M'$, if there is a chain of
elementary Nielsen moves which transforms $M$ to $M'$.
\end{defn}

It is easy to see that if $M\sim_N M'$ then $M$ and $M'$ generate
the same subgroup of $G$. For this reason Nielsen equivalence is
a very useful tool for studying the subgroup structure of various
groups.

We prove the following statement which can be regarded as an
``acylindrical accessibility'' result for finitely generated
groups acting on real trees. Indeed, our theorem says that there
is a bound on the size of a ``fundamental domain'' for a minimal
$D$-acylindrical isometric action of a $k$-generated group on an
$\mathbb R$-tree:

\begin{thm}\label{D} Let $G$ be a freely indecomposable finitely
generated group acting by isometries on an $\mathbb R$-tree $X$.
Let $D\ge 0$. Suppose that $G\ne 1$ is  not infinite  cyclic, that
the action of $G$ is $D$-acylindrical, nontrivial (does not have a
fixed point) and minimal (that is $X$ has no proper $G$-invariant
subtrees).

Let $\varepsilon>0$ be an arbitrary real number. Then any finite
generating $k$-tuple $Y$ of $G$ is Nielsen-equivalent to a
$k$-tuple $S$ such that:

\begin{enumerate}
\item There is a finite subtree $T_{\varepsilon}\subseteq X$ of
measure at most $2D(k-1)+\varepsilon$ such that
$GT_{\varepsilon}=X$;
\item for some $x\in T_{\varepsilon}$ we
have
\[
\max\{d(x,sx)| s\in S\} \le 2D(k-1)+\varepsilon.
\]
\end{enumerate}

\end{thm}

By the \emph{measure of $T_{\varepsilon}$} we mean the sum of the lengths of intervals
in any subdivision of $T_{\varepsilon}$ as a disjoint union of finitely many
intervals. This is equal to the 1-dimensional Hausdorff measure
of $Y$. If $X$ is a simplicial tree and $T_{\varepsilon}$ is a simplicial subtree,
then the measure of $T_{\varepsilon}$ is the number of edges in $T_{\varepsilon}$.

Theorem~\ref{D} immediately implies the following since actions
with trivial arc stabilizers are $0$-acylindrical:

\begin{cor}\label{E}
Let $G$ be a finitely generated freely indecomposable non-cyclic
group which acts by isometries on an $\mathbb R$-tree $X$ with
trivial arc stabilizers. Then for any $\varepsilon>0$ and any
finite generating tuple $Y$ of $G$ there is a tuple $S$
Nielsen-equivalent to $Y$ and a point $x\in X$ such that
$d(x,sx)\le \varepsilon$ for all $s\in S$.
\end{cor}

Thus if the action of $G$ is $0$-acylindrical, that is arc
stabilizers are trivial, then a finite generating set of $G$ can
be made by Nielsen transformations to have arbitrarily small
translation length.   Not surprisingly our methods let us recover
the same bound $c(G,D)=2D(k-1)$ on the complexity of acylindrical
accessibility splittings as the one given in \cite{W}. The main
ingredient is a theory of Nielsen methods for groups acting on
hyperbolic spaces that we systematically developed in
\cite{KW,KW1}. This theory is analogous to Weidmann's treatment of
actions on simplicial trees~\cite{W}, but the case of arbitrary
hyperbolic spaces is technically much more complicated. Note that
the proof of Theorem~\ref{D} completely avoids the Rips machinery
for groups acting on $\mathbb R$-trees~\cite{BF95,Ka01} and
Theorem~\ref{D} makes no traditional stability assumptions about
the action. Rather, we use the fact that an $\mathbb R$-tree is
$\delta$-hyperbolic for any $\delta>0$, which allows us to make a
limiting argument for $\delta$ tending to zero.

We thank the referee whose detailed and insightful comments have helped to greatly improve this paper.

\section{The main technical tool}

Our main tool is a technical result (see Theorem~\ref{mainA}
below) obtained by Kapovich and Weidmann in \cite{KW1}. It is
motivated by the Kurosh subgroup theorem (see \cite{LS77, Bass93})
for free products, which states that a subgroup of a free product
$\ast_{i\in I} A_i$ is itself a free product of a free group and
subgroups that are conjugate to subgroups of the factors $A_i$.

\begin{defn}
Suppose that $(X,d)$ is an $\mathbb R$-tree and that $U$ is a finitely generated group that acts on
$X$ by isometries. Put $E(U):=\{x\in X\,|\, ux=x\hbox{ for some
}u\in U-1\}$. If $U$ does not fix a point of $X$, let $X_U$ be the
minimal $U$-invariant subtree of $X$. If $U$ fixes a point of $X$,
let $X_U$ denote the set of all points of $X$ fixed by $U$.

Finally, put $X(U)$ to be the smallest $U$-invariant subtree
containing $X_U\cup E(U)$. Thus $X(U)$ is a nonempty $U$-invariant
subtree of $X$.
\end{defn}

We generalize the notion of Nielsen equivalence as follows. The
objects which correspond to the tuples of elements of $G$ are the
{\em $G$-tuples:}

\begin{defn}[$G$-tuple]
Let $G$ be a group.

Let $n\ge 0$, $m\ge 0$ be integers such that $m+n>0$. We will say
that a tuple $M=(U_1,\dots, U_n;H)$ is a {\em $G$-tuple} if $U_i$ is
a non-trivial subgroup of $G$ for $i\in\{1,\ldots ,n\}$ and
$H=(h_1,\ldots ,h_m)\in G^m$ is an $m$-tuple of elements of $G$.
We will denote $\overline{M}=U_1\cup \dots \cup U_n\cup
\{h_1,\ldots ,h_m\}$ and call $\overline{M}$ the \emph{underlying
set} of $M$. Note that $\overline M$ is nonempty since $m+n>0$.
\end{defn}

By analogy with the Kurosh subgroup theorem we will sometimes
refer to the subgroups $U_i$ as {\em elliptic components} of $M$.
This is justified since in most applications of our methods, in
particular the proof of Theorem~\ref{D}, the subgroups $U_i$ are
generated by sets of elements with short translation length. We
should stress, however, that $U_i$ need not be fixing a point of a
tree on which $G$ acts. We will also refer to $H$ as the {\em
hyperbolic component} of $M$. We have the following notion of
equivalence for $G$-tuples which generalizes the classical Nielsen
equivalence.

\begin{defn}[Equivalence of $G$-tuples]
We will say that two $G$-tuples $M=(U_1,\dots, U_n;H)$ and
$M'=(U_1',\dots , U_n'; H')$ are \emph{equivalent} if
$H=(h_1,\ldots ,h_m)$ and $H'=(h_1',\ldots ,h_m')$ and $M'$ can be
obtained from $M$ by a chain of moves of the following type:
\begin{enumerate}
\item For some $1\le j\le n$ replace $U_j$ by  $gU_jg^{-1}$ where
\[
g\in\langle \{h_1,\ldots ,h_m\}\cup U_1\cup\ldots \cup
U_{j-1}\cup U_{j+1}\cup\ldots \cup U_n\rangle .
\]

\item For some $1\le i\le n$ replace $h_i$ by
$h_i'=g_1h_ig_2$ where
\[
g_1,g_2\in \langle \{h_1,\ldots ,h_{i-1},h_{i+1},\ldots
,h_m\}\cup U_1\cup\ldots\cup U_n\rangle .
\]

\end{enumerate}
\end{defn}

Our main technical tool is Theorem~\ref{mainA} stated below. It
is a corollary of Theorem~2.4 in~\cite{KW1} that deals with
arbitrary group actions on $\delta$-hyperbolic geodesic metric
spaces. 

\begin{thm}[Kapovich-Weidmann]\cite{KW1}\label{mainA}
Let $G$ be a group acting on an $\mathbb R$-tree $(X,d)$ by
isometries.

Let  $M=(U_1,\dots, U_n;H)$ be a $G$-tuple where $H=(h_1,\ldots
,h_m)$  and let
\[U=\langle \overline{M}\rangle=\langle
\{h_1,\ldots ,h_m\}\cup \bigcup_{i=1}^n U_i\rangle \le G.\] Let
$\varepsilon>0$ be an arbitrary positive number.

Then either $U=U_1\ast \dots U_n\ast F(H)$ or there exists a
$G$-tuple $M'=(U_1',\dots, U_n';H')$ with $H'=(h_1',\ldots ,h_m')$
such that $M'$ is equivalent to $M$ and at least one of the
following holds:

\begin{enumerate}
\item $d(X(U_i'),X(U_j'))< \varepsilon$ for some $1\le i<j\le n$.

\item $d(X(U_i'),h_j'X(U_i'))< \varepsilon$ for some
$i\in\{1,\ldots ,n\}$ and $j\in\{1,\ldots ,m\}$.

\item There exists a point $x\in X$ such that $d(x,h_jx)<
\varepsilon$ for some $j\in\{1,\ldots ,m\}$.
\end{enumerate}
\end{thm}

\begin{rem}
The statement of Theorem~2.4 in \cite{KW1} does not use the tree
$X(U)$ defined above but, rather, the tree $X_{\delta}(U)$ (where
$\delta\ge 0$) defined as follows.

Let
\[
E_{\delta}(U):=\{x\in X: d(ux,u)\le 100\delta \text{ for some }
u\in U, u\ne 1.\}
\]

Then the tree $X_{\delta}(U)$ is defined as $E_{\delta}(U)$ if $U$
fixes a point and as the smallest $U$-invariant subtree of $X$
containing $X_U$ and $E_{\delta}(U)$ otherwise.

Note that $X(U)\subseteq X_{\delta}(U)$ by construction. Let $x\in
E_{\delta}(U)$ and let $y\in X(U)$ be such that
$d(x,y)=d(x,X(U))$. Let $u\in U, u\ne 1$ be such that $d(x,ux)\le
100\delta$.  Note that the choice of $y$ guarantees that there exists no $z\in [x,y]$ with $y\neq z$ and $uz=z$ as  otherwise $z\in X(U)$ and $d(x,z)<d(x,y)$. It follows that $d(x,ux)=2d(x,y)+d(y,uy)$ and hence
$d(x,y)\le 50\delta$. Thus we have shown that $E_{\delta}(U)$ is contained in the $50\delta$-neighborhood of $X(U)$.
This implies that, whether or not $U$ fixes a point, 
 for any $\delta\ge 0$ the trees $X(U)$
and $X_{\delta}(U)$ are $50\delta$-Hausdorff close.

Since an $\mathbb R$-tree is $\delta$-hyperbolic for any
$\delta>0$, Theorem~2.4 of \cite{KW1} now directly implies
Theorem~\ref{mainA} above by taking the limit $\delta\to 0$.
\end{rem}

  We deploy Theorem~\ref{mainA} in the proof of Theorem~\ref{D}
  for a "generator transfer"
process to analyze a freely indecomposable subgroup generated by a
finite set $Y=\{y_1,\ldots ,y_k\}$ with $k$ elements. First we
start with a $G$-tuple $M_1=(;H_Y)$ where $H_Y=(y_1,\ldots ,y_k)$.
We then construct a sequence of $G$-tuples $M_1,M_2,\dots $ by
repeatedly applying Theorem~\ref{mainA} in order to either "drag"
elements of the "hyperbolic" components of $G$-tuples into their
"elliptic" components or to join two elliptic components to form
one new elliptic component. A simple observation shows that the
length of the sequence $M_1,M_2,\dots $ is bounded by $2k-1$. The
desired result is then obtained by analyzing the terminal member
of this sequence.

\section{Groups acting on real trees}

 We define "generating trees" exactly as in \cite{BW}.

\begin{defn}[Generating tree]
Let $U=\langle S\rangle$ be a group acting by isometries on an
$\mathbb R$-tree $X$. We say that a tree $T_U\subset X$ is an
$S$-{\em generating tree of $U$} if $Y\cap sY\neq\emptyset $ for
all $s\in S$. We further say that $T_U$ is a \emph{generating tree
of $U$} if $T_U$  is an $S$-generating tree for $U$ for some
generating  set $S$ of $U$.
\end{defn}

The following lemma is an immediate consequence of the
definitions; it also justifies the term ``generating tree''.

\begin{lem} Let $X$ be an $\mathbb R$-tree and
let $U=\langle S\rangle$ be a group acting on $X$ by isometries.
Let $T_U\subset X$ be an $S$-generating tree of $U$. Then the following hold:

\begin{enumerate}
\item The set $UT_U$ is connected and $U$-invariant, and hence is
a subtree of $X$.

\item  If $U$ does not fix a point of $X$, then $UT_U$ contains
the minimal $U$-invariant subtree $X_U$ of $X$.
\end{enumerate}
\end{lem}

\medskip We can now observe that for acylindrical actions the fix
point set $E(U)$ cannot be too far from the minimal $U$-invariant tree
$X_U$, i.e. that $X(U)$ and $X_U$ are close. The following is a simple
exercise.

\begin{lem}\label{t2} Let $U$ be a group acting by isometries on an $\mathbb
R$-tree $X$ and suppose this action is $D$-acylindrical for some
$D\ge 0$.   Then the following hold:.

\begin{enumerate}

\item Let $y\in X$ be such that $d(y,X_U)=R$. Then for any $u\in
U, u\ne 1$ we have $d(y,uy)\ge 2R-2D$.

\item The set $E(U)$ is contained in the $D$-neighborhood of
$X_U$. Moreover, $X(U)$ and $X_U$ are $D$-Hausdorff close.

\item If $U$ fixes a point $x\in X$ then for any $y\in X(U)$ we have
  $d(x,y)\le D$.

\end{enumerate}

\end{lem}

\medskip Before we proceed with the proof of Theorem~\ref{D} let us recall some
more notions from \cite{KW1}.

\begin{defn}[Partitioned tuple]
Let $G$ be a group.  If $Y$ is an $n$-tuple of elements of
$G$ we will say that $n$ is \emph{the length of} $Y$ which we
denote by $L(Y)$. We will say that $M=(Y_1,\dots, Y_p; H)$ is a
\emph{partitioned tuple} in $G$ if $p\ge 0$ and $Y_1,\dots, Y_p,
H$ are finite tuples of elements of $G$ such that at least one of
these tuples has positive length and such that for any $i\ge 1$
the tuple $Y_i$ has positive length.

We will call the sum of the lengths of $L(Y_1)+\dots + L(Y_p)+L(H)$
\emph{the length of} $M$ and denote it by $L(M)$. We further call the
$L(M)$-tuple, obtained by concatenating the tuples $Y_1,\ldots,Y_p,H$,
\emph{the tuple underlying $M$}.
\end{defn}

\begin{defn}[Complexity]
  Let $M=(Y_1,\dots, Y_n; H)$ be a partitioned tuple.  As in
  \cite{W,KW1}, we define the \emph{complexity} of $M$ to be the pair
  $(L(H),n)$.  Thus the complexity is an element of ${\mathbb N}^2$
  where $\mathbb N=0,1,2,\ldots$.  We order $\mathbb N^2$ by setting
  $(m,n)\le (m',n')$ if $m<m$ or if $m=m'$ and $n\le n'$. This gives a
  well-ordering on $\mathbb N^2$ and allows us to compare
  complexities.
\end{defn}

\begin{rem}\label{fixtuples} To any partitioned tuple $M=(Y_1,\dots,
  Y_p; H)$ we associate the $G$-tuple $\tau=(U_1,\ldots ,U_p;H)$,
  where $U_i=\langle Y_i\rangle$.  Suppose now that $\tau$ is
  equivalent to a $G$-tuple $\tau'=(U_1',\ldots , U_p',H')$. The
  definition of equivalence of $G$-tuples implies that there is a
  partitioned tuple $M'=(Y_1',\dots, Y_p'; H')$ with associated
  $G$-tuple $\tau'$ such that the tuples underlying $M$ and $M'$ are
  Nielsen-equivalent. Moreover, we can choose $Y_i'$ to be conjugate
  to~$Y_i$ for each $1\le i\le p$.
\end{rem}

\medskip\noindent{\em Proof of Theorem \ref{D}.}

Recall that in Theorem~\ref{D} $Y$ is a given $k$-tuple generating
$G$. Clearly, it is enough to prove the statement of Theorem~\ref{D}
under the assumption that $Y$ is not Nielsen-equivalent to a tuple
containing $1\in G$. Thus we will assume that every $k$-tuple
Nielsen-equivalent to $Y$ consists of nontrivial elements.

Now in order to prove Theorem~\ref{D} it suffices to establish:

{\bf Claim.} There exist a $k$-tuple $S$ Nielsen-equivalent to $Y$ and
an $S$-generating tree $T_{\varepsilon}$ of measure at most
$2D(k-1)+\varepsilon$.

Let $N=(S_1,\ldots S_n;H)$ be a partitioned tuple of elements of $G$.
Let $k_i=L(S_i)$ for $1\le i\le n$. We say that $N$ is {\em good} if
$U_i=\langle S_i\rangle\ne 1$ for all $i\ge 1$ and if for each $U_i,
i\ge 1$ there exists an $S_i$-generating tree $T_i$ of measure at most
$2D(k_i-1)+\frac{2k_i-1}{2k}\varepsilon$.

We define $N_1$ to be the partitioned tuple $N_1=(;Y)$. Clearly $N_1$
is good.

\medskip

\noindent{\bf Choice of $M$}. Let $M=(S_1,\ldots S_n;H)$ with
$H=(h_1,\ldots ,h_m)$ be a partitioned tuple of minimal complexity
among all partitioned good tuples with the underlying tuple being
Nielsen equivalent to $Y$. The partitioned tuple $N_1$ satisfies the
above qualifying constraints and hence such an $M$ exists.  \medskip

We will show that $M=(S_1;-)$. This would immediately imply the Claim.

Suppose that $M$ is not of this type. Recall that $G$ is freely
indecomposable and not infinite cyclic. It follows from
Theorem~\ref{mainA} and Remark~\ref{fixtuples} that there exists a
good partitioned tuple $M'=(S_1',\ldots S_n',H')$ of the same
complexity as $M$ with $H'=(h_1',\ldots ,h_m')$ such that the
underlying tuple of $M'$ is Nielsen equivalent to $Y$ and such that
the following holds. If we denote $U_i'=\langle S_i'\rangle$ for $1\le
i\le n$ then at least one of the following occurs:

\begin{enumerate}
\item $d(X(U_1'),X(U_2'))\le \frac{\varepsilon}{2k}$;
  
\item $d(X(U_1'),h_m'X(U_1'))\le\frac{\varepsilon}{2k}$;
  
\item $d(y,h_m'y)\le \frac{\varepsilon}{2k}$ for some $y\in X$.
\end{enumerate}

Recall that  $k_i=L(S_i)=L(S_i')$ for $1\le i\le n$.

\medskip {\bf Case 1.} Suppose that $d(X(U_1'),X(U_2'))\le
\frac{\varepsilon}{2k}$. Choose $x_1\in X(U_1')$ and $x_2\in X(U_2')$
such that $d(x_1,x_2)\le \frac{\varepsilon}{2k}$. By part (2) of
Lemma~\ref{t2} there is $y_i\in X_{U_i'}$ such that $d(y_i,x_i)\le D$
for $i=1,2$. It follows that $d(y_1,y_2)\le
2D+\frac{\varepsilon}{2k}$. Since by assumption $M'$ is good, we can
choose an $S_i'$-generating tree $T_{U_i'}$ for $U_i'$ of measure at
most $2D(k_i-1)+\frac{2k_i-1}{2k}\varepsilon$.

Let $i\in\{1,2\}$. If $U_i'$ does not fix a point then
$X_{U_i'}\subset U_i'T_{U_i'}$ and there exists a $u_i\in U_i'$ such
that $y_i\in u_iT_{U_i'}$. If $U_i'$ fixes a point, then
$X_{U_i'}$ is the fixed set of $U_i'$ and hence $y_i\in X_{U_i'}$ is
fixed by $U_i'$. In this case we can assume that in $M'$ we have
$T_{U_i'}=\{y_i\}$. This is clearly an $S_i'$-generating tree for
$U_i'$ of measure zero. With $u_i=1$ we also still have $y_i\in
u_iT_{U_i'}$.

Denote $S_i'':=u_iS_i'u_i^{-1}$ for $i=1,2$. Then $\langle
S_i''\rangle =\langle S_i'\rangle =U_i'$ and $u_iT_{U_i'}$ is the
$S_i''$-generating tree for $U_i'$.

Moreover $S_i''$ is Nielsen-equivalent to $S_i'$ since $u_i\in U_i'$
for $i=1,2$. Put $V:=\langle U_1', U_2'\rangle$.  Then
$T_{V}=u_1T_{U_1'}\cup [y_1,y_2]\cup u_2T_{U_2'}$ is a generating tree
for $V$ with respect to $S_1''\cup S_2''$. The measure of $T_{V}$ is
at most
\begin{gather*}
  2D(k_1-1)+\frac{2k_1-1}{2k}\varepsilon+2D(k_2-1)+\frac{2k_2-1}{2k}\varepsilon+2D+
  \frac{\varepsilon}{2k}=\\
  2D(k_1+k_2-1)+\frac{2(k_1+k_2)-1}{2k}\varepsilon.
\end{gather*}

Hence the partitioned tuple $M'':=(S_1''\cup S_2'',S_3',\ldots ,S_n';
H')$ is good. Since the underlying tuple of $M''$ is Nielsen
equivalent to $Y$ and since $M''$ has smaller complexity than $M'$ and
$M$, we obtain a contradiction to the choice of $M$.

\medskip {\bf Case 2.}  Suppose that $d(X(U_1'), h_m'X(U_1'))\le
\frac{\varepsilon}{2k}$.  As in (1) we see that there exist $y_1\in
X_{U_1'}$ and $y_2\in h_m'X_{U_1'}$ such that $d(y_1,y_2)\le
2D+\frac{\varepsilon}{2k}$.

{\bf Subcase 2A.} Suppose first that $U_1'$ does not fix a point.
Again, since $M'$ is good, after replacing $S_1'$ by a conjugate (in
$U_1'$) tuple we can assume that there exists an $S_1'$-generating tree
$T_{U_1'}$ for $U_1'$ of measure at most
$2D(k_1-1)+\frac{2k_1-1}{2k}\varepsilon$ such that $y_1\in T_{U_1'}$.
Clearly $(h_m')^{-1}y_2$ lies in $X_{U_1'}$. Since $X_{U_1'}\subseteq
U_1'T_{U_1'}$ it follows that $u_2'(h_m')^{-1}y_2\in T_{U_1'}$ for
some $u_2'\in U_1'$.

Hence $T_{V}=T_{U_1'}\cup [y_1,y_2]$ is a generating tree for the
subgroup $V=\langle S_1',u_2' (h_m')^{-1}\rangle=\langle S_1',
h_m'\rangle$ with respect to the generating set $S_1''= S_1'\cup\{u_2'
(h_m')^{-1}\}$.

Moreover, the measure of $T_{V}$ is at most

\begin{align*}
  &2D(k_1-1)+\frac{2k_1-1}{2k}\varepsilon
  +2D+\frac{\varepsilon}{2k}\le
  2D((k_1+1)-1)+\frac{2(k_1+1)-1}{2k}\varepsilon.
\end{align*}

This implies that the tuple $M''=(S_1'',S_2',\ldots ,S_n';H'')$ is
good, where $H''=(h_1',\ldots, h_{m-1}')$. Since $M''$ has smaller
complexity than does $M$, this contradicts the choice of $M$.

{\bf Subcase 2B.}  Suppose now that $U_1'$ fixes a point.  Then the
assumption $d(X(U_1'), h_m'X(U_1'))\le \frac{\varepsilon}{2k}$ implies
that there exist $y_1\in X(U_1')$ and $y_2\in h_m'X(U_1')$ such that
$d(y_1,y_2)\le \frac{\varepsilon}{2k}$.  Note that since $U_1'$ fixes
a point, $X_{U_1'}$ is the fixed set of $U_1'$. Denote
$y_2'=(h_m')^{-1}y_2\in X(U_1')$. Let $x$ be a point fixed by $U_1'$.
By Lemma~\ref{t2} we have $d(y_1,x)\le D$ and $d(y_2',x)\le D$.  Put
\[
K=[y_2',x]\cup [x,y_1]\cup [y_1,y_2]
\]
 
Then $x$ is fixed by $U_1'$ and hence by $S_1'$, while $h_m'y_2'=y_2$.
Thus $K$ is an $S_1''$-generating tree for the subgroup $V=\langle
S_1', h_m'\rangle$ where $S_1''=(S_1',h_m')$. Note that
$L(S_1'')=L(S_i')+1\ge 2$ by construction. The tree $K$ has measure at
most
\[
2D+\frac{\varepsilon}{2k}\le 2D(2-1)+\frac{2\cdot 2-1}{2k}\varepsilon
\]

Again put $H''=(h_1',\ldots, h_{m-1}')$ and $M''=(S_1'',S_2',\ldots
,S_n';H'')$. Then $M''$ is good and has smaller complexity than $M$.
This contradicts the choice of $M$.

\medskip {\bf Case 3.} Suppose that $d(y,h_m' y)\le
\frac{\varepsilon}{2k}$ for some $y\in X$.  In this case we replace
$M'$ by
\[
M''=(S_1',\ldots ,S_n', S_{n+1}';H'')
\]
where $S_{n+1}'=(h_m')$ and $H''=(h_1',\ldots ,h_{m-1}')$.

By assumption on $Y$ we have $h_m'\ne 1$. Hence $M''$ is good and of
smaller complexity than $M$, which again yields a contradiction.
\hfill$\Box$

\begin{rem}
  Essentially the same argument as in the proof of Theorem~\ref{D}
  implies a relative version of our main result. If $U_1,\dots, U_n$
  are subgroups of a group $G$, we will say that $G$ is \emph{freely
    indecomposable relative to $U_1,\dots, U_n$} if there does not
  exist an action of $G$ on a simplicial tree $X$ with trivial edge
  stabilizers and without inversions such that each $U_i$ fixes a
  vertex. Essentially the same argument as in the proof of
  Theorem~\ref{D} implies the following relative version of our main
  result.
  
  Suppose $G$ is freely indecomposable relative to $U_1, .., U_n$ and
  that $G$ is generated by $U_1\cup\dots \cup U_n
  \cup\{s_1,...,s_k\}$. Suppose $G$ acts $D$-acylindrically on an
  $\mathbb R$-tree in such a way that each $U_i$ fixes a point. Then
  for any $\epsilon>0$ there exists a $G$-generating tree of measure
  at most $2D(n+k-1)+\epsilon$.
  
\end{rem}

 \end{document}